\documentclass[a4paper,10pt]{article}
\usepackage[utf8x]{inputenc}
\usepackage{amssymb}

\title{Lagrangian geometry of algebraic varieties}
\author{ Nikolay Tyurin \\
BLTPh JINR (Dubna) and Volga Region  Mat Center (Kazan)}

\date{}

\begin{document}

\maketitle

\begin{abstract} Every algebraic variety can be regarded as a symplectic manifold being equipped
with a Kahler form. Therefore it is natural to study lagrangian geometry of any algebraic variety.
We present two basic constructions which can be applied to a sufficiently wide set of algebraic varieties.
\end{abstract}

{\bf Lagrangian geometry.} Every compact algebraic variety $X$ by the very definition (see [1]) admits principal polarization
so a very ample line bundle $L$ which defines an embedding $X  \hookrightarrow \mathbb{C} \mathbb{P}^N$
to a projective space. Therefore if  a standard Kahler form on the last one is fixed it induces the corresponding
Kahler form on $X$; the form is not unique but the corresponding lagrangian geometry essentially depends on the
cohomology class $c_1(L) \in H^2(X, \mathbb{Z})$ only.

 Thus the choice $c_1(L)$ on $X$ discovers a new direction in the studies of $X$: taking the corresponding Kahler form
$\omega$ as a symplectic form one can pose the following questions. First, one says that a real submanifold $S \subset X$
is lagrangian if $\omega|_S$ identically vanishes and $\rm{dim}_{\mathbb{C}} X = \rm{dim}_{\mathbb{R}} S$. And the
first classification question arises as follows: which homology classes $\nu_i \in H_n(X, \mathbb{Z})$
are realizable by smooth lagrangian submanifolds? And the next basic question is what are possible topological types
of smooth lagrangian submanifolds?

{\bf Example:} the complex projective plane $\mathbb{C} \mathbb{P}^2$ admits essentially unique principal polarization given by
 $L = {\cal O}(1)$ which corresponds to the standard Kahler form for the Fubini - Study metric. Just recently the answer to the
  basic question of the lagrangian geometry has been completed for this case:  one proves that the Klein bottle can not
  appear as a smooth lagrangian submanifold, and the complete list of possible lagrangian submanifolds contains real projective
  plane $\mathbb{R} \mathbb{P}^2$ and torus $T^2$ both representing the trivial homology class.

 More specific questions arise when one starts the classification programme for lagrangian submanifolds: if two lagrangian submanifolds
  $S_1$ and $S_2$ can be joined by a family of lagrangian submanifolds $S_t$ then one says that $S_1$ and $S_2$ are equivalent up
  to lagrangian deformation; if it exists a Hamiltonian function $H(x, t)$ such that the corresponding flow $\phi^t_{X_H}$
  moves $S_1$ to $S_2$ then these are called Hamiltonian isotopic. The last equivalence is much more delicate then the first one.
  Then a high level classification question for a given algebraic variety with a fixed principal polarization asks how many
  equivalence classes (up to lagrangian deformation or up to Hamiltonian isotopy) exists for a given class $\nu_i \in H_n(X, \mathbb{Z})$
  and a given topological type $\rm{top} S$. It is not hard to see that the answer does not depend on the particular choice of the Kahler form
  but on the class $c_1(L)$ only.

  {\bf Example.} Consider the case of lagrangian tori in $\mathbb{C} \mathbb{P}^2$. Since the projective plane is toric then
  it exists the family of Liouville tori which are lagrangian by the Liouville theorem; at the same time one of them
  is the famous Clifford torus. Untill 1996 one thought that no other lagrangian tori of the same periods but Hamiltonian non isotopic
  do not exist, and Yu. Chekanov construction (see [2]) presented an example of lagrangian torus which has the same periods but
  which is not Hamiltonian isotopic to the Clifford torus. This torus was called exotic.

{\bf Generalization of Chekanov construction.}  More general, every {\it toric} variety by the very definition admits Liouville type tori: such a variety is presented by
  the phase space of completely integrable system, and the common level sets for the first integrals are lagrangian tori. Thus
  toric geometry is a subject where algebraic geometry meets symplectic geometry: T. Delzant proved that all the data
  are encoded by a simple combinatoric object, a convex polytop $P_X$, and both the symplectic geometry and algebraic
  geometry can be recovered from $P_X$ (see [3]). The Chekanov construction can be translated to the algebraic language as follows:
  the components of the boundary divisor (= preimage of the boundary $\partial P_X$) can be combined in the elements of a pencil
  of algebraic divisors such that each element is invariant under the Hamiltonian action of one first integral, say, $f_1$.
  Therefore for $\mathbb{C} \mathbb{P}^2$ instead of a pair of first integrals $(f_1, f_2)$ one considers the pair $(f_1, \{Q\})$,
  where the second element is a pencil of conics on $\mathbb{C} \mathbb{P}^2$: the integrability requires that the Hamiltonian action
  of $f_1$ preserves the fibers of the pencil. Then the base of the pencil $\mathbb{C} \mathbb{P}^1$ admits two points $p_{\pm}$ which underly
  singular fibers, and every smooth loop $\gamma \subset \mathbb{C} \mathbb{P}^1 \backslash \{p_{\pm}\}$ gives a lagrangian torus
  in whole $\mathbb{C} \mathbb{P}^2$; if $\gamma$ is non contractible then the resulting torus is of the Clifford type,
  if $\gamma$ is contractible then the torus is exotic (or of the Chekanov type).

   Therefore we can generalize the Chekanov  construction for any toric algebraic variety $X$. First of all we need to extract
  certain information from the convex polytop $P_X \subset \mathbb{R}^n$. It is formed by $r$ linear inequalities
  bounded by facets $\Delta_i$ of $P_X$ such that the Picard group $\rm{Pic} X$ subjects the exact sequence
  $$
  0 \to \mathbb{Z}^n \to \mathbb{Z}^r \to \rm{Pic} X \to 0,
  $$
  where the second  $\mathbb{Z}^r$ is spanned on the components $D_i$ of the boundary divisor $D_b$. Then each point
  of $\mathbb{R}^n$ gives a relation on $D_i$ of the form
  $$
  \sum_{i =1}^r a_i [D_i] = 0 \quad \rm{in} \quad \rm{Pic} X.
  $$
  Since each $D_i$ is effective, realized by a subvariety, all the coefficients  $a_i$ can not be of the same signs, and we can rearrange the last equality as follows
  $$
  \sum_{a_i \geq 0} a_i [D_i] = \sum_{a_i < 0} \vert a_i \vert [D_i],
  $$
  denoting the corresponding divisors as $D_+$ and $D_-$; by the very construction they represent the same class $[D_+ ] = [D_-]$
  in the Picard group. Therefore if we take the corresponding line bundle $L \to X$ such that $c_1(L) = \rm{P.D.} [D_{\pm}]$
  then it exists a pencil $<D_+, D_-> \subset \vert L \vert$ in the complete linear system, defined by holomorphic sections of $L$.
  Excluding the base set $B = D_+ \cap D_- $ from $X$ one gets a map $\psi: X \backslash B \to \mathbb{C} \mathbb{P}^1$;
  it is not hard to see that the points $p_{\pm}= \psi(D^0_{\pm})$ exhaust the set of points  underly singular fibers.
  At the same time the choice of the pencil $<D_+, D_->$ derives a linear condition on the complete set of first integrals $(f_1, ..., f_n)$
  such that it exists a set $(\tilde f_1, ..., \tilde f_{n-1})$ of integrals such that each element of the pencil is invariant
  with respect to the Hamiltonian action of each $\tilde f_i$. 

   Therefore for a toric algebraic variety $X$ one can fix data $(\tilde f_1, ..., \tilde f_{n-1}, \psi)$ choosing a direction
  in $\mathbb{R}^n$; for a fixed data every smooth loop $\gamma$ on the complement $ <D_+, D_-> \backslash \{p_{\pm} \}$ and
  a generic value set $c_1, ..., c_{n-1}$ for functions $\tilde f_i$ we have the union
  $$
  T_{\gamma, c_i} = \cup_{p \in \gamma} (\psi^{-1}(p) \cap N(c_1, ..., c_{n-1})),
  $$
  where $N(c_1, ..., c_{n-1})$ is the common level set for the functions $\tilde f_i$. The {\bf Main theorem} says that  $T_{\gamma, c_i}$
  is a smooth lagrangian torus: if $\gamma$ is non contractible in $<D_+, D_-> \backslash \{ p_{\pm} \}$ then the torus is
  of the standard type, if $\gamma$ is contractible then the resulting torus is of exotic type, see [4]. Moreover,
  in [5] one proves that if $X$ admits a standard monotone lagrangian torus then it exists an exotic lagrangian torus
  which is monotone as well.
  
   At the same time one can apply the construction to a non toric algebraic variety: in [4] one presents the structure $(f_1, ..., f_{n-1}, \psi)$
  on flag variety $F^3$, which is the full flag in $\mathbb{C}^3$. This variety is non toric, however it admits natural non complete set
  of first integrals $(f_1, f_2)$ and it is possible to find a pencil $\psi: F^3 \backslash B \to \mathbb{C} \mathbb{P}^1$ such that
  the fibers are invariant under the Hamiltonian action of $f_i$. In this non toric situation the base $\mathbb{C} \mathbb{P}^1$ contains three points which underly
  singular fibers therefore we have more complicated picture since a smooth loop on the complement "base minus three points"
  can represent different classes being non contractible. In this setup one can reconstruct the famous Gelfand - Zeytlin lagrangian sphere,
  see [4]; different types of smooth lagrangian tori can be derived from the picture, see [5].
  
  {\bf Generalization of Mironov construction.} In [6] A. Mironov presented new example of (Hamiltonian) minimal lagrangian submanifolds in $\mathbb{C}^n$ and
  $\mathbb{C} \mathbb{P}^n$. The geometric essence of the Mironov construction is the following: one uses together the standard toric structures and the natural real structures  
which exist on both the varieties. From this point of view it is not hard to find a generalization of the Mironov construction: it can be applied to an algebraic
variety which admits a real structure and a (possible non complete) toric action which is transversal to the real structure. If $X$ is our algebraic variety
and $\omega$ is compatible with the real structure then the real part $X_{\mathbb{R}} \subset X$ must be isotropical with respect to $\omega$. If moreover
the real part has maximal possible dimension $\rm{dim}_{\mathbb{R}}  X_{\mathbb{R}} = \rm{dim}_{\mathbb{C}} X$ then $X_{\mathbb{R}}$ is lagrangian.

Suppose we have a Hamiltonian toric action spanned by moment maps (or first integrals) $f_1, ..., f_k$. Then for a generic value set $c_1, ..., c_k$
one takes the common level set $N(c_1, ..., c_k) = \{ f_i = c_i \}$ and the intersection $S_{\mathbb{R}}(c_1, ..., c_k) = X_{\mathbb{R}} \cap
N(c_1, ..., c_k)$. Since the toric action is transversal to the real structure by the assumption the last  submanifold $S_{\mathbb{R}}(c_1, ..., c_k)$
is transversal to the Hamiltonian vector fields $X_{f_i}$, and if we apply the toric action, generated by $X_{f_i}$, to the last submanifold
we get $T^k(S_{\mathbb{R}}(c_1, ..., c_n)) \subset X$ which is a subcycle (since possibly it has self intersections).
In [7] one proves the {\bf Main theorem} which states that $T^k(S_{\mathbb{R}}(c_1, ..., c_k))$ is a lagrangian immersion
(in the smooth case -- lagrangian submanifold) in $X$. Using different sets of $f_i$'s one presents different types of lagrangian submanifolds
in the Grassmannian $\rm{Gr}(2, 4)$. At the end in [7] one proposes that the Grassmannian $\rm{Gr}(2, n)$ admits at least $n+1$ topologically
distinct lagrangian submanifolds.

Indeed, the Grassmann variety $\rm{Gr}(k, n)$ gives a natural example of algebraic variety which admits natural real structure and toric structure of the desired types.
As the space of $k-1$ - dimensional projective subspaces of $\mathbb{C} \mathbb{P}^{n-1}$ it inherits the toric action of $T^{n-1}$ coming from the standard
toric action on $\mathbb{C} \mathbb{P}^{n-1}$; on the other hand it admits the natural real structure such that the real part $\rm{Gr}_{\mathbb{R}}(k, n)$
has right dimension. For the construction one can take any subtorus in $T^{n-1}$, which leads to the natural grading by the rank of this subtorus,
starting with $X_{\mathbb{R}}$ itself, which is claimed to be of the zero level.  In [8] one realizes the programme for the level 1 case.  {\bf Main theorem} here reads as follows:
 Mironov cycle is smooth Lagrangian submanifold, isomorphic to topologically non trivial
 fiber bundle over real Grassmannian ${\rm Gr}_{\mathbb{R}}(k, n -1)$,  where the fiber is either $S^1 \times S^{k-1}$ for even $k$ or generalized Klein bottle for
  for odd  $k$. In particular  for $\rm{Gr}(2, n)$  for any $n$ it gives certain topologically non trivial $T^2$ - bundle over ${\rm Gr}_{\mathbb{R}}(2, n-1)$.

  {\bf References:}
  
  [1] P. Griffits, J. Harris, {\it "Principles of algebraic geometry"}, NY, Wiley, 1978;
  
  [2] Yu. Chekanov, {\it "Lagrangian tori in a symplectic vector space and global symplectomorphisms"}, Math. Z., 223: 4 (1996), 547–559;
  
  [3] T. Delzant, {\it "Hamiltoniens périodiques et images convexes de l’application moment"}, Bull. Soc. Math. France, 116: 3 (1988), 315–339;
  
  [4] N. Tyurin, {\it "Pseudotoric structures: Lagrangian submanifolds and Lagrangian fibrations"}, Russian Mathematical Surveys (2017), 72: 3, 513 - 546;
  
  [5] N. Tyurin, {\it "Monotonic Lagrangian Tori of Standard and Nonstandard Types in Toric and Pseudotoric Fano Varieties"}, Proc. Steklov Inst. Math., 307 (2019), 267–280;
  
  [6] A. Mironov, {\it  "New examples of Hamilton-minimal and minimal Lagrangian manifolds in $\mathbb{C}^n$ and $\mathbb{C} \mathbb{P}^n$"}, Sb. Math., 195: 1 (2004), 85–96;
  
  [7] N. Tyurin, {\it "Mironov Lagrangian cycles in algebraic varieties"}, Sb. Math., 212:3 (2021), 389–398; 
  
  [8] N. Tyurin, {\it "Examples of Mironov cycles in Grassmann varieties"}, Siberian Math. J., 62: 2 (2021),  457–465.

\end{document}